\documentclass{amsart}
\theoremstyle{plain}
\newtheorem{theorem}{Theorem}[section]
\newtheorem{proposition}[theorem]{Proposition}
\newtheorem{corollary}[theorem]{Corollary}
\newtheorem{lemma}[theorem]{Lemma}

\newtheorem{definition}[theorem]{Definition}

\newtheorem{remark}[theorem]{Remark}

\newcommand{\bfC}{{\mathbb C}}

\newcommand{\bfR}{{\mathbb R}}

\newcommand{\bari}{{\overline i}}
\newcommand{\barj}{{\overline j}}
\newcommand{\bark}{{\overline k}}
\newcommand{\barl}{{\overline \ell}}
\newcommand{\barpartial}{{\overline \partial}}

\newcommand{\barz}{{\overline z}}

\newcommand{\barbeta}{{\overline \beta}}

\newcommand{\tildeg}{{\widetilde g}}
\newcommand{\tildeu}{{\widetilde u}}

\newcommand{\mapright}[1]{\smash{\mathop{   \hbox to 0.7cm{\rightarrowfill}}
  \limits^{#1}}}

\begin{document}

\title{\Large Harmonic total Chern forms and stability}
\author{Akito Futaki}
\address{Department of Mathematics, Tokyo Institute of Technology, 2-12-1,
O-okayama, Meguro, Tokyo 152-8551, Japan}
\email{futaki@math.titech.ac.jp}
\date{March 24, 2006 }
\begin{abstract} In this paper we will perturb the scalar curvature of compact K\"ahler manifolds
by incorporating it with higher Chern forms, and then show that the perturbed
scalar curvature has many common properties with the unperturbed scalar curvature.
In particular the  perturbed
scalar curvature becomes a moment map, with respect to a perturbed symplectic structure, 
on the space of all complex structures on a
fixed symplectic manifold, which extends the results of Donaldson and Fujiki
on the unperturbed case.
\end{abstract}
\maketitle

\section{Introduction}

Many works have been done on the relationship between the existence of constant
scalar curvature K\"ahler metrics and stability in the sense of geometric invariant
theory. 
A way of seeing this relationship is through the moment map picture of an infinite
dimensional set up as done by Donaldson \cite{donaldson97} and Fujiki \cite{fujiki92}.
They showed that the set of all K\"ahler metrics with constant scalar curvature becomes 
the zero set of the moment map for the action of
the group of Hamiltonian symplectomorphisms on the space of all compatible complex structures on
a fixed symplectic manifold. Recall that for a Hamiltonian action of a compact Lie
group $K$ on a compact K\"ahler manifold, having a zero of the moment map along an orbit of the 
complexified group $K^c$-action is equivalent to the stability of the orbit of the reductive group $K^c$
(c.f. \cite{donkro}, section 6.5). Applying this fact in finite dimensions to the infinite 
dimensional space of all compatible complex structures we see a relationship between the existence of 
constant scalar curvature
K\"ahler metrics and infinite dimensional symplectic-GIT stability.

The purpose of this paper is to perturb the scalar curvature by incorporating it with higher Chern
classes, and show that the perturbed scalar curvature shares many common properties with the
unperturbed scalar curvature. Especially the set of all K\"ahler metrics with constant 
perturbed scalar curvature is the zero set of the moment map with respect to a perturbed
symplectic form on the space of all compatible complex structures on a fixed symplectic manifold.
This extends the earlier results of Donaldson and Fujiki in the unperturbed case.

Let $M$ be a compact symplectic manifold with a fixed symplectic form $\omega$ and of 
dimension $2m$.
Let $\mathcal{J}$ be the set of all $\omega$-compatible integrable complex structures.
Then for each $J \in {\mathcal J}$, $(M, \omega, J)$ becomes a K\"ahler manifold. For a pair $(J,t)$ of a complex
structure $J$ and a small real number $t$, define a smooth function $S(J,t)$ on $M$ by

\begin{equation} \frac{S(J,t)}{2m\pi}\ \omega^m = c_1(J)\wedge \omega^{m-1} + t c_2(J) \wedge \omega^{m-2} +
\cdots + t^{m-1} c_m(J) \end{equation}
where $c_i(J)$ is the $i$-th Chern form with respect to
the K\"ahler structure $(\omega, J)$ on $M$, i.e. they are defined by
\begin{equation}  \det ( I + \frac i{2\pi}t\Theta) = 1 + tc_1(J) + \cdots + t^mc_m(J), \end{equation}
$\Theta$ being the curvature matrix of the Levi-Civita connection. Note that $S(J,0)$ is equal to the trace of the Ricci curvature
$g^{i\barj}R_{i\barj}$ which is one half of the Riemannian scalar curvature. But since $S(J,0)$ 
more often appears in the computations in K\"ahler geometry than the Riemannian scalar curvature does, 
we will call $S(J,0)$ the scalar curvature in this paper. 
We also call $S(J,t)$ {\sl the perturbed scalar curvature}. As mentioned above 
the main result of this paper is to show that the perturbed scalar curvature becomes a
moment map on $\mathcal J$ with respect to some symplectic structure (Theorem 2.2 in
the next section).  

This paper is organized as follows. In section 2, we will prove Theorem 2.2. 
We will give two proofs along the lines of \cite{donaldson97} and \cite{tianbook00}.
In section 3, we study the analogy to extremal K\"ahler metrics in our perturbed case.
We will see that the perturbed extremal K\"ahler metrics are critical points of the functional
on ${\mathcal J}$ given by the squared $L^2$-norm of the perturbed scalar curvature but not critical
points of the functional on the space of K\"ahler forms given by the same integral. In section 4 we will
recall Bando's result \cite{bando83} on the obstructions to the existence of K\"ahler metrics with
harmonic higher Chern classes and study the relevant Mabuchi functional in the perturbed case.
In section 5, we will give a deformation theory of extremal K\"ahler metrics to the perturbed extremal
K\"ahler metrics extending earlier results of LeBrun and Simanca \cite{lebrunsimanca93}, \cite{lebrunsimanca94}.

\section{Perturbed symplectic structure on the space of complex structures}

Let $(M, \omega)$ be a compact symplectic manifold of dimension $2m$ and $\mathcal J$ 
the space of all $\omega$-compatible complex structures on $M$. This means that $J \in \mathcal J$
if and only if $ \omega(JX, JY) = \omega(X, Y)$ for all vector fields $X$ and $Y$, and 
$\omega(X, JX) > 0$ for all non-zero $X$. 
For later purposes it is convenient to assume that $J$ acts on the cotangent bundle
rather than the tangent bundle. Fixing $J \in \mathcal J$,
we decompose the complexified cotangent bundle into holomorphic and anti-holomorphic parts,
i.e. $\pm \sqrt{-1}$-eigenspaces of $J$:
\begin{equation}
T^{\ast}M\otimes\bfC = T_J^{\ast\prime}M \oplus T_J^{\ast\prime\prime}M, 
\qquad T_J^{\ast\prime\prime}M = \overline{T_J^{\ast\prime}M}.
\end{equation}
Taking arbitrary $J^{\prime} \in \mathcal J$ we also have the decomposition with respect to $J^{\prime}$
\begin{equation}
T^{\ast}M\otimes\bfC = T_{J^{\prime}}^{\ast\prime}M \oplus T_{J^{\prime}}^{\ast\prime\prime}M, 
\qquad T_{J^{\prime}}^{\ast\prime\prime}M = \overline{T_{J^{\prime}}^{\ast\prime}M}.
\end{equation}
If $J^{\prime}$ is sufficiently close to $J$ then $T_{J^{\prime}}^{\ast\prime}M $ 
can be expressed as a graph over $ T_J^{\ast\prime}M$ as
\begin{equation}
T_{J^{\prime}}^{\ast\prime}M = \{\ \alpha + \mu(\alpha) \ |\ \alpha \in T_J^{\ast\prime}M\ \}
\end{equation}
for some endomorphism $\mu$ of $T_J^{\ast\prime}M$ into $T_J^{\ast\prime\prime}M$:
\begin{eqnarray}
 \mu &\in& \Gamma (\mathrm{End}(T_J^{\ast\prime}M, T_J^{\ast\prime\prime}M)) \nonumber \\
 &\cong&  \Gamma (T_J^{\prime}M \otimes T_J^{\ast\prime\prime}M) \cong \Gamma(T_J^{\prime}M\otimes T_J^{\prime}M)
\end{eqnarray}
where in the last identification we used the K\"ahler metric defined by the pair $(\omega, J)$. This can be
expressed in the notation of tensor calculus with indices as 
$$ \mu^i{}_{\bark} \mapsto g^{j\bark}\mu^i{}_{\bark}=: \mu^{ij} $$
where we chose a local holomorphic coordinate system
$(z^1, \cdots, z^m)$ and wrote $\omega$ as $\omega = \sqrt{-1}\ g_{i\barj} dz^i \wedge d\overline{z^j}$.

\begin{lemma} With the above identification understood, $\mu$ lies in the symmetric part 
$\Gamma(\mathrm{Sym}(T_J^{\prime}M\otimes T_J^{\prime}M))$ of $\Gamma(T_J^{\prime}M\otimes T_J^{\prime}M)$.
\end{lemma}

\begin{proof} The symplectic form $\omega$ gives a natural identification between the tangent bundle and 
the cotangent bundle. This identification then gives a natural symplectic structure on the cotangent
bundle, which we denote by $\omega^{-1}$. If $\omega$ is $J$-invariant, then $\omega^{-1}$ is also
$J$-invariant. 
For the complex structure $J$, $\omega^{-1}$ is expressed in terms of the K\"ahler metric of 
the K\"ahler structure $(\omega, J)$ as
$$\omega^{-1} = - \sqrt{-1}\ g^{i\barj}\frac{\partial}{\partial z^i} \wedge \frac{\partial}{\partial\overline{z^j}},$$
where we used the local expression of $\omega$ as above.
Since $\omega^{-1}$ is $J$-invariant and any 1-forms $\alpha$ and $\beta$ in $T_J^{\ast\prime}M$ are eigenvectors
of $J$ belonging to $\sqrt{-1}$, we have
$$ \omega^{-1}(\alpha,\beta) = 0. $$
Similarly we have
$$ \omega^{-1}(\mu\alpha,\mu\beta) = 0 $$
and, since $\omega^{-1}$ is also $J^{\prime}$-invariant, we also have
$$ \omega^{-1}(\alpha + \mu\alpha,\ \beta + \mu\beta ) = 0. $$
Thus we obtain
\begin{equation}
\omega^{-1}(\alpha, \mu\beta) = \omega^{-1}(\beta, \mu\alpha)
\end{equation}
which implies that $\mu \in \Gamma(T_J^{\prime}M\otimes T_J^{\prime}M)$ is symmetric
because in the local expression, 
\begin{equation}
\mu^{ji}\alpha_i\beta_j = \mu^{ij}\alpha_i\beta_j,
\end{equation}
as desired.
\end{proof}

Considered infinitesimally, 
the tangent space $T_J\mathcal J$ to 
$\mathcal J$ at $J$ is a subspace of $\mathrm{Sym}(T_J^{\prime}M\otimes T_J^{\prime}M)$.

Then the $L^2$-inner product on 
$\mathrm{Sym}(T_J^{\prime}M\otimes T_J^{\prime}M)$ gives $\mathcal J$ 
a K\"ahler structure. But we perturb this K\"ahler structure in the following way.
Let $t$ be a small real number. For $\mu$ and $\nu$ in the tangent space $T_J\mathcal J$, we define
\begin{equation}
(\nu, \mu)_t = \int_M mc_m(\overline{\nu}_{jk}
\,\mu^i{}_\barl \frac {\sqrt{-1}}{2\pi}\,
dz^k \wedge d\overline{z^{\ell}}, \omega\otimes I + \frac {\sqrt{-1}}{2\pi}\,t\Theta,
\cdots, \ \omega\otimes I + \frac {\sqrt{-1}}{2\pi}\,t\Theta)
\end{equation}
where $c_m$ is the polarization of the determinant viewed as a $GL(m,\bfC)$-invariant 
polynomial, i.e. $\ c_m(A_1, \cdots, A_m)$ is the coefficient of $m!\, t_1 \cdots t_m$ in $\det(t_1A_1 + \cdots +
t_mA_m)$, where $I$ denotes the identity matrix and $\Theta = \barpartial(g^{-1}\partial g)$ 
is the curvature form of the
Levi-Civita connection, and where 
$u_{jk}\mu^i_{\bar l}$  should be understood as
      the endomorphism of  $T'_JM$ which sends 
      $\partial/\partial z^j$  to  $u_{jk}\mu^i_{\bar l}{\partial/\partial z^i}$.

Note that
$$ c_m(A,\ \cdots\ ,A) = \det A. $$
This is similar to the wedge product
$$ \alpha_1\wedge\ \cdots \wedge\alpha_m$$
for the type $(1,1)$-forms $ \alpha_1,\ \cdots , \alpha_m$.
For we have 
$$ \alpha\ \wedge\ \cdots\ \wedge\alpha = \det (a_{i\barj})dz^1\wedge d\barz^1
\wedge \cdots \wedge dz^m \wedge d\barz^m$$
when  $ \alpha = \sum a_{i\barj}dz^i\wedge d\barz^j$.
Therefore there is a symmetry between
the endomorphism part and the form part in the integration of (9).
This symmetry will be used in this work and was used in the work of Bando \cite{bando83} quoted in
the next section.

When $t=0$, $(\cdot,\cdot)_t$ gives the standard $L^2$-inner product which is anti-linear in the first factor $\nu$ and
linear in the second factor $\mu$. If the real number $t$ is sufficiently small, $(\cdot,\cdot)_t$ is still
positive definite. 

Let $\mathcal G$ be the group of all Hamiltonian symplectomorphisms of $(M, \omega)$. The
Lie algebra of $\mathcal G$ is isomorphic to the Poison algebra $C^{\infty}_0(M)$ 
of all smooth functions on $M$ with average $0$:
$$ C^{\infty}_0(M) = \{\ u \in C^{\infty}(M)\ |\ \int_M\ u\ \omega^m = 0\}. $$
$\mathcal G$ acts on $\mathcal J$ as holomorphic isometries. 

\begin{theorem}\ \ For each fixed small real number $t$, $S(J,t)/2m\pi$ gives 
an equivariant moment map
on $\mathcal J$ if we consider $S(J,t)/2m\pi$ as an element of the dual space of  $C^{\infty}_0(M)$ 
by the pairing
$$ <\frac{S(J,t)}{2m\pi}, u> = \int_M u\ \frac{S(J,t)}{2m\pi}\ \omega^m. $$
\end{theorem}

The case $t=0$ is due to Donaldson (\cite{donaldson97}) and Fujiki (\cite{fujiki92}), and 
a mildly different proof in this case was also given in Tian's book \cite{tianbook00}. 

To prove the theorem, let us consider two operators
$$ P : C^{\infty}_0 \to T_J {\mathcal J}, $$
$$ Q: T_J {\mathcal J} \to C^{\infty}_0(M),$$
where $P$ represents the infinitesimal action of the Lie algebra $C^{\infty}_0$ on 
${\mathcal J}$ via Hamiltonian action and $Q$ represents the derivative of the map
which associates to $J \in {\mathcal J}$  
the perturbed scalar curvature $\frac 1{2m\pi} S(J,t)$ of the K\"ahler manifold
$(M, \omega, J)$. We need to show
$$ \Re (P(u), \sqrt{-1}\mu)_t = <Q(\mu),u> $$

To compute $P(u)$, we have only to compute $L_X J$ for a smooth vector field $X$.

\begin{lemma} For a smooth vector field $X = X^{\prime} + X^{\prime\prime}$ we have
$$ L_X J = 2\sqrt{-1} \nabla_J^{\prime\prime} X^{\prime}
 - 2\sqrt{-1} \nabla_J^{\prime} X^{\prime\prime}. $$
In particular, if $X_u$ is the Hamiltonian vector field of $u$,
 $$ P(u) = 2\sqrt{-1} \nabla_J^{\prime\prime} X^{\prime}_u.$$
\end{lemma}

\begin{proof}
Since $(L_X J)\alpha = L_X(J\alpha) - JL_X \alpha$, if $\alpha$ is a type $(1,0)$-form,
\begin{equation}
(L_X J)\alpha = \sqrt{-1}(L_X \alpha - (L_X \alpha)^{1,0} + (L_X\alpha)^{0,1}) 
= 2\sqrt{-1} (L_X \alpha)^{0,1}.
\end{equation}
On the other hand
\begin{equation}
L_X\alpha = d(\alpha(X^{\prime})) + i(X)(\partial_J\alpha + \barpartial_J\alpha).
\end{equation}
Thus 
\begin{equation}
(L_X\alpha)^{0,1} = \barpartial_J (\alpha(X^{\prime})) 
+ i(X^{\prime})\barpartial_J\alpha.
\end{equation}
But
\begin{eqnarray*}
\barpartial_J (\alpha(X^{\prime}))  &=& \nabla_J^{\prime\prime}(\alpha(X^{\prime})) \\
&=& (\nabla_J^{\prime\prime} \alpha)(X^{\prime}) + \alpha(\nabla_J^{\prime\prime} X^{\prime}) 
= (\barpartial_J \alpha)(X^{\prime}) + \alpha(\nabla_J^{\prime\prime} X^{\prime}).
\end{eqnarray*}
This implies  
\begin{equation}
\barpartial_J (\alpha(X^{\prime})) + i(X^{\prime})(\barpartial_J \alpha) = 
\alpha(\nabla_J^{\prime\prime} X^{\prime}) 
\end{equation}

From (10), (12) and (13) we get
\begin{equation}
(L_XJ) \alpha = \alpha (2\sqrt{-1}\  \nabla_J^{\prime\prime}X^{\prime}).
\end{equation}
Similarly, if $\alpha$ is a $(0,1)$-form, then
\begin{equation}
(L_XJ) \alpha = \alpha (-2\sqrt{-1}\  \nabla_J^{\prime}X^{\prime\prime}).
\end{equation}
From (14) and (15) we get the lemma. This completes the proof.
\end{proof}

From this lemma we get for the real function $u$

\begin{eqnarray}
&&\Re (P(u), \sqrt{-1}\mu)_t = 2 \Re (\nabla_J^{\prime\prime}X^{\prime}_u, \mu)_t
 \\
&=& 2 \Re \int_M mc_m(u_{jk}
\,\mu^i{}_\barl \frac  {\sqrt{-1}}{2\pi}
dz^k \wedge d\overline{z^{\ell}}, \omega\otimes I + \frac {\sqrt{-1}}{2\pi}\,t\Theta,
\cdots, \ \omega\otimes I + \frac  {\sqrt{-1}}{2\pi}\,t\Theta). \nonumber
\end{eqnarray}

Next we need to compute $Q$. We will do this in two ways along the lines of \cite{donaldson97} and 
\cite{tianbook00}. First we follow the arguments of \cite{donaldson97} 
just word for word.

If identify $T_{J^{\prime}}^{\ast\prime}M$ with 
$T_J^{\ast\prime}M$ through $\alpha + \mu\alpha \mapsto \alpha$,
this identification induces identifications of differential forms with all degrees, 
which we denote by 
$ \iota : \Omega_{J^{\prime}}^{p,q} \to \Omega_J^{p,q}$.

\begin{lemma}
With the above identification we have the following.
\begin{enumerate}
\item
If a 1-form $\gamma = \alpha + \barbeta \in 
T_J^{\ast\prime}M  \oplus T_J^{\ast\prime\prime}M$ 
is written also as 
$\gamma = \alpha^{\prime} + \mu\alpha^{\prime} + \overline{\beta^{\prime} + 
\mu\beta^{\prime}} \in 
T_{J^{\prime}}^{\ast\prime}M  \oplus T_{J^{\prime}}^{\ast\prime\prime}M$ 
then 
$$\overline{\beta^{\prime}} = \barbeta - \mu\alpha$$ 
up to first order in $\mu$. Namely 
$$\iota(\overline{\beta^{\prime} + \mu\beta^{\prime}}) = \barbeta - \mu\alpha$$
up to first order in $\mu$. 
\item If a fixed 2-form $\chi = \chi^{2,0} + \chi^{1,1} + \chi^{0,2} \in 
\Omega_J^{2,0} \oplus  \Omega_J^{1,1} \oplus \Omega_J^{0,2}$ 
has $\chi^{\prime 1,1}$ as a $(1,1)$-component with respect to $J^{\prime}$,
then 
$$ \iota(\chi^{\prime 1,1}) = \chi^{1,1} - \mu\chi^{2,0} - \mu\chi^{0,2} $$
up to first order in $\mu$, 
where we extended the operation of $\mu$ to higher degree tensors in the
obvious way. 
\end{enumerate}
\end{lemma}

Hereafter we use the notation $\equiv$ to mean "up to first order in $\mu$".

\begin{proof}\ \ (a)\ \ From $\alpha^{\prime} = \alpha - \overline{\mu\beta^{\prime}}$
we see 
\begin{equation*}
\overline{\beta^{\prime}} = \barbeta - \mu\alpha^{\prime} 
= \barbeta - \mu(\alpha - \overline{\mu\beta^{\prime}}) 
\equiv \barbeta - \mu\alpha. 
\end{equation*}
(b)\ \ 
If a fixed 2-form is written
also as $\chi = (1+\mu)\alpha_1\wedge (1+\mu)\alpha_2 + 
(1+\mu)\alpha_3 \wedge \overline{(1+\mu)\beta_1} +
\overline{(1+\mu)\beta_2} \wedge \overline{(1+\mu)\beta_3}
 \in 
\Omega_{J^{\prime}}^{2,0} \oplus  \Omega_{J^{\prime}}^{1,1} \oplus 
\Omega_{J^{\prime}}^{0,2}$, then a similar computation as in the proof of (a) shows
\begin{eqnarray*}
\alpha_3 \wedge \barbeta_1 
&\equiv& \chi^{1,1} - \alpha_1 \wedge \mu\alpha_2 - \mu\alpha_1 \wedge \alpha_2 
- \overline{\mu\beta_1 \wedge\beta_2} - \overline{\beta_1 \wedge \mu\beta_2}\\
&\equiv& \chi^{1,1} - \mu\chi^{2,0} - \mu\chi^{0,2}.
\end{eqnarray*}
This completes the proof.
\end{proof}

\begin{corollary}\ \ Let $E \to M$ be a vector bundle. If $\nabla$ is a fixed connection
of $E$ and $\nabla = \nabla_J^{\prime} + \nabla_J^{\prime\prime}$ with respect to
the complex structure $J$, then by the identification above
$\nabla_{J^{\prime}}^{\prime\prime}$ is identified
with $\nabla_J^{\prime\prime} - \mu\nabla_J^{\prime}$ up to first order in $\mu$.
\end{corollary}

\noindent
{\it Proof of Theorem 2.2}\ \ The identification $\iota : T_{J^{\prime}}^{\ast\prime}M \to
T_J^{\ast\prime}M$ is a Hermitian isometry up to first order in $\mu$, and we can 
consider the Levi-Civita connections $\nabla_J$ and $\nabla_{J^{\prime}}$ as
two unitary connections on the same bundle. If $J$ is fixed and $\nabla^{\prime\prime}$
is varied by $\sigma \in \Omega^{0,1}(\mathrm{End}(T^{\prime}M))$ then the connection
changes by $\sigma - \sigma^{\ast}$. On the other hand, if a connection 
$\nabla = \nabla_J^{\prime} +\nabla_J^{\prime\prime}$ is fixed and $J$ varies 
to $J^{\prime}$ by $\mu$, then the new $\nabla_{J^{\prime}}^{\prime\prime}$ is
identified with $\nabla_J^{\prime\prime} - \mu\nabla_J^{\prime}$ up to first order
in $\mu$ by Corollary 2.5.

 Now we compute $\nabla_{J^{\prime}}^{\prime\prime}$ for a 1-form $\alpha$
 of $T_J^{\ast\prime}M$, which is strictly speaking equal to
$\iota\circ\nabla_{J^{\prime}}^{\prime\prime}\circ \iota^{-1}(\alpha)$.
But $\nabla_{J^{\prime}}^{\prime\prime}\circ \iota^{-1}(\alpha)$ is 
$\Omega^{0,1}_{J^{\prime}}$-part of $d(\alpha + \mu\alpha)$ up to first order
in $\mu$.
From this and 
Lemma 2.4, (b), we get
\begin{equation}
\nabla_{J^{\prime}}^{\prime\prime}\alpha \equiv \nabla_J^{\prime\prime} \alpha
+ \nabla_J^{\prime}(\mu\alpha) -\mu(\nabla_J^{\prime}\alpha).
\end{equation}
On $T_J^{\ast\prime}M \otimes T_J^{\ast\prime}M$, $\mu$ acts as a derivation.
To make the notations clear we will denote by $\mu_1$ (resp. $\mu_2$) the
action of $\mu$ on the first (resp. second) factor. So, on 
$T_J^{\ast\prime}M \otimes T_J^{\ast\prime}M$, we have $\mu = \mu_1 \otimes 1 +
1\otimes \mu_2$. With these notations the right hand side of (17) is equal to
\begin{eqnarray}
\nabla_J^{\prime\prime}\alpha + \mu_2\nabla_J^{\prime}\alpha + (\nabla_J^{\prime}\mu)\alpha
- \mu(\nabla^{\prime}\alpha) 
&\equiv& \nabla_J^{\prime\prime}\alpha - \mu_1\nabla_J^{\prime}\alpha + (\nabla_J^{\prime}\mu)\alpha \\
&\equiv& (\nabla_J^{\prime\prime} - \mu\nabla_J^{\prime})\alpha + (\nabla_J^{\prime}\mu)\alpha.
\nonumber
\end{eqnarray}

By Corollary 2.5, $\nabla_J^{\prime\prime} - \mu\nabla_J^{\prime}$ is the expression
under our identification of 
$J^{\prime}$-$(0,1)$-component of a fixed connection $\nabla_J$. Thus the variation of
the Levi-Civita connection is $\sigma - \sigma^{\ast}$ where $\sigma = \nabla_J^{\prime}\mu$.
Notice that $\sigma$ must be a $(0,1)$-form with values in $\mathrm{End}(T_J^{\prime}M)$.
So, in local expressions
$$ \nabla_J^{\prime}\mu = (\nabla_j\mu^i{}_{\barl}d\overline{z^{\ell}})$$
with $i$ column index, $j$ row index. 
Since it is convenient to distinguish the covariant derivative as the endomorphism part 
from the covariant exterior derivative as the form part, we shall write $\nabla_J$ 
to denote the covariant derivative as the endomorphism part and $d^{\nabla_J}$
to denote the covariant exterior derivative as the form part. 
Thus, under the variation $\delta J = \mu$ of the complex structure,
the variation $\delta\Theta$ of the curvature
matrix $\Theta$ is
$$ \delta\Theta = d^{\nabla_J}(\sigma - \sigma^{\ast}). $$
Its $(1,1)$-part is
$$ (\delta\Theta)^{1,1} =d^{ \nabla_J^{\prime}}(\nabla_J^{\prime} \mu) - 
(d^{\nabla_J^{\prime}}(\nabla_J^{\prime} \mu))^{\ast}. $$
\noindent
Since the exterior covariant derivative 
$d^{\nabla_J}(\omega\otimes I\ + \ \frac{\sqrt{-1}}{2\pi}\ t\Theta)$ of 
$\omega\otimes I\ + \ \frac{\sqrt{-1}}{2\pi}\ t\Theta$ vanishes, we have

\begin{eqnarray*}
&&\delta \int_M u\ \frac{S(J,t)}{2m\pi} \omega^m\\
&& =2\Re \int_M u\ mc_m(\frac{\sqrt{-1}}{2\pi}d^{\nabla_J^{\prime}}(\nabla_J^{\prime}\mu),
\omega\otimes I + \frac{\sqrt{-1}}{2\pi}\ t\Theta,\ \cdots, \omega\otimes I + 
\frac{\sqrt{-1}}{2\pi}\ t\Theta) \\
&& = -2\Re \int_M mc_m(\frac{\sqrt{-1}}{2\pi}\overline{d^{\nabla_J^{\prime\prime}}u}\wedge 
\nabla_J^{\prime}\mu, 
\omega\otimes I + \frac{\sqrt{-1}}{2\pi}\ t\Theta,\ \cdots, \omega\otimes I + 
\frac{\sqrt{-1}}{2\pi}\ t\Theta) 
\end{eqnarray*}
Now the invariant polynomial $c_m$ takes determinant for the endomorphism part, and therefore
we may interchange the roles of the form part and the endomorphism part in the integration above.
Thus by the vanishing of $d^{\nabla_J}(\omega\otimes I\ + \ \frac{\sqrt{-1}}{2\pi}\ t\Theta)$
again we can use integration by parts for the covariant derivative of the endomorphism part.
Hence we have 
$$
\delta \int_M u\ \frac{S(J,t)}{2m\pi} \omega^m = 
2\Re \int_M mc_m(\frac{\sqrt{-1}}{2\pi}\overline{\nabla_J^{\prime\prime}d^{\nabla_J^{\prime\prime}}u}\wedge \mu,
\omega\otimes I + \frac{\sqrt{-1}}{2\pi}\ t\Theta,\ \cdots, \omega\otimes I 
+ \frac{\sqrt{-1}}{2\pi}\ t\Theta).
$$
where the term 
$\frac{\sqrt{-1}}{2\pi}\overline{\nabla_J^{\prime\prime}d^{\nabla_J^{\prime\prime}}u}\wedge \mu$ is
expressed in local coordinates
$$ \frac{\sqrt{-1}}{2\pi}u_{kj} dz^k \wedge \mu^i{}_{\barl}dz^{\barl},$$
where $u_{kj} = \nabla_j \nabla_k u$.
This coincides with (16), completing the proof of Theorem 2.2.

\medskip

\noindent
{\it Alternate proof of Theorem 2.2}\ \ We only need to show that $<Q(\mu),u>$ is equal
to (16). To compute $Q$ we take a local coordinates $(x^1, \cdots,\ x^{2m})$ with respect
to which $\omega$ is the standard symplectic form on $\bfR^{2m}$, by using Darboux's
theorem. Let $J_t$ be a family of complex structures with $J_0 = J$. Then we have
$$\dot{J}|_{t=0} = 
2\sqrt{-1}\mu - 2\sqrt{-1}\overline{\mu}. $$
This follows because, by taking the derivative of 
$$J_t (\alpha + \mu(t)\alpha) = \sqrt{-1}(\alpha + \mu(t)\alpha)
$$ 
with $\dot{\mu}(0) = \mu$, we have
$$ \dot{J}(\alpha) = 2\sqrt{-1}\mu. $$
Let $g_t = \omega J_t$ be the Riemannian metric induced by $J_t$. Then the
Christoffel symbols of $g_t$ are written as
$$ \Gamma_{t,jk}^i = \frac12 g_t^{i\ell}\left(
\frac{\partial g_{t,\ell j}}{\partial x^k} + \frac{\partial g_{t,\ell k}}{\partial x^j}
- \frac{\partial g_{t,j k}}{\partial x^{\ell}} \right). $$
At $p \in M$ we may assume that $g_{ij}(p) = \delta_{ij},\ dg_{ij}(p) = 0$, and 
$$ J(p) = \left(\begin{array}{cc} O & -I \\ I & O \end{array}\right)$$
where $g = g_0$. Then $ \Gamma_{t,jk}^i$ is of order $t$, and 

\begin{eqnarray*}
R_{t, ijk\ell} &=& g_t( 
\nabla_{\frac{\partial}{\partial x^i}}\nabla_{\frac{\partial}{\partial x^j}}\frac{\partial}{\partial x^{\ell}}
- 
\nabla_{\frac{\partial}{\partial x^j}}\nabla_{\frac{\partial}{\partial x^i}}\frac{\partial}{\partial x^{\ell}},
\frac{\partial}{\partial x^k}) \\
&=& g_{t,sk}\frac12 \left(\frac{\partial^2 g_{t,pj}}{\partial x^i \partial x^{\ell}} - 
\frac{\partial^2 g_{t,j{\ell}}}{\partial x^i \partial x^p}
 - \frac{\partial^2 g_{t,pi}}{\partial x^j \partial x^{\ell}} + 
\frac{\partial^2 g_{t,i{\ell}}}{\partial x^j \partial x^p} \right)g_t^{ps}\\
&& +\ \ 
\mathrm{quadratic\ terms\ in\ the\ first\ derivatives\  of}\ g.
\end{eqnarray*}
Taking the derivative with respect to $t$ at $t = 0$,
\begin{eqnarray*}
\frac d{dt} |_{t=0} R_{t, ijk\ell}  &=& \frac12 (\dot{g}_{kj,i\ell} 
- \dot{g}_{j\ell,ik} - \dot{g}_{ki,j\ell} + \dot{g}_{i\ell,jk}).
\end{eqnarray*}

Now we compute the right hand side in terms of local holomorphic coordinates $z^1,\ \cdots, z^m$.
The only terms involved in the integration are $\dot{g}_{\bari\barl,jk}$'s and their complex
conjugates, and we also have
$$ \dot{g}_{\bari\barl} = - \sqrt{-1}g_{\bari p}2\sqrt{-1}\mu^p_{\barl} = 2\mu_{\bari\barl}.$$
Thus
$$
 \frac12 \dot{g}_{\bari \barl, jk} \sqrt{-1}dz^k \wedge d\overline{z^{\ell}}
= \mu_{\bari \barl, jk} \sqrt{-1}dz^k \wedge d\overline{z^{\ell}}.
$$
Hence we get 
$$ <Q(\mu),u> = 2\Re \int_M u\ mc_m(\frac{\sqrt{-1}}{2\pi}
d^{\nabla_J^{\prime}}(\nabla_J^{\prime}\mu),\ 
\omega\otimes I + \frac{\sqrt{-1}}{2\pi}\ t\Theta,\ \cdots,
\ \omega\otimes I + \frac{\sqrt{-1}}{2\pi}\ t\Theta). $$
As in the last part of the previous proof this last term coincides with (16). 
This completes the alternate proof.\\

\section{Perturbed extremal K\"ahler metrics}

For a real or complex valued smooth function $u$ on a K\"ahler manifold $(M,g)$ we put
$$ \mathrm{grad}^{\prime} u = \sum_{i,j = 1}^m 
g^{i\barj}\frac{\partial u}{\partial z^{\barj}}
\frac{\partial}{\partial z^i}$$
and call it the gradient vector field of $u$. Strictly speaking the real part of 
$ \mathrm{grad}^{\prime} u$ is the gradient vector field of $u$, but we identify
a real vector field with its $T^{\prime}M$-part.

\begin{definition}\ \ A K\"ahler metric $g = (g_{i{\barj}})$ is
said to be a perturbed extremal K\"ahler metric 
if the gradient vector field 
$$
{\mathrm grad}^{\prime} S(J,t)  = \sum_{i,j = 1}^m 
g^{i\barj}\frac{\partial S(J,t)}{\partial z^{\barj}}
\frac{\partial}{\partial z^i}
$$
of the perturbed scalar curvature $S(J,t)$ is a holomorphic vector field. 
\end{definition}

\begin{proposition}\ \ Critical points of the functional
$$J \mapsto \int_M S(J,t)^2 \omega^m$$
on $\mathcal J$ are perturbed extremal K\"ahler metrics.
\end{proposition}

\begin{proof}\ \ Let $J(s)$ be a smooth family of complex structures such that
$J(0) = J$ and $\dot{J}(0) = \mu$. By the proof of Theorem 2.2
$$\left. \frac{d}{ds}\right|_{s=0}  \int_M u\ S(J(s),t)\ \omega^m =
2m\pi \Re (\nabla^{\prime\prime}\nabla^{\prime\prime}u,
 \mu)_t$$
 for all real smooth function $u$ with $\int_M u\,\omega^m = 0$.
We take $u$ to be $v:= S(J,t) - \int_M S(J,t) \omega^m/\int_M \omega^m$
and $\mu$ to be $(-\sqrt{-1})$-times the infinitesimal action of the Hamiltonian vector field of
$v$ at $J$. Then using the 
above equality and Lemma 2.3
$$\left.\frac{d}{ds}\right|_{s=0}  \int_M v\ S(J(s),t) = 2m\pi \Re (\nabla^{\prime\prime}\nabla^{\prime\prime}u, \mu)_t.$$
From this we get
\begin{eqnarray*}
\left.\frac{d}{ds}\right|_{s=0}  \int_M S(J(s),t)^2\, \omega^m &=& 
2 \int_M S(J,t) \left.\frac{d}{ds}\right|_{s=0} S(J(s),t)\, \omega^m \\
&=& 2 \int_M v \left.\frac{d}{ds}\right|_{s=0} S(J(s),t)\, \omega^m \\
&=& 4m\pi \Re  (\nabla^{\prime\prime}\nabla^{\prime\prime}u, \mu)_t.
\end{eqnarray*}
This shows that $J$ is a critical point if and only if 
$$\nabla^{\prime\prime}\mathrm{grad}^{\prime}S(J,t) = 0,$$
i.e. the K\"ahler metric of $(M, \omega, J)$ is a perturbed extremal K\"ahler
metric.
\end{proof}

\begin{remark}\ \ In the case of unperturbed extremal K\"ahler metrics when $t=0$, 
such K\"ahler metrics are also the critical points of the functional
$$ \omega \mapsto \int_M S(\omega)^2 \omega^m $$
on the space of all K\"ahler forms $\omega$ in a fixed K\"ahler class $[\omega_0]$
where $S(\omega)$ denotes the scalar curvature of the K\"ahler form $\omega$,
(c.f. \cite{calabi85}).
But when $t \ne 0$
the perturbed extremal K\"ahler metrics are not the critical points of the functional
$$ \omega \mapsto \int_M S(\omega,t)^2 \omega^m$$
on the space of all K\"ahler forms in a fixed K\"ahler class where 
\begin{eqnarray}
\frac{S(\omega,t)}{2m\pi}\,\omega^m &=& c_1(\omega)\wedge \omega^{m-1} + 
t\,c_2(\omega)\wedge \omega^{m-2} + \cdots + t^{m-1}c_m(\omega)\\
&=& \frac 1t (\det(\omega\otimes I + \frac {\sqrt{-1}}{2\pi}\,t\Theta) - \omega^m),\nonumber
\end{eqnarray}
$c_j(\omega)$ being the $j$-th Chern form with respect to $\omega$:
$$ \det (1 + t\ \frac{\sqrt{-1}}{2\pi}\Theta) = 1 + tc_1(\omega) + \cdots + t^{m-1}c_m(\omega).$$
Note that
we use the notation $S(\omega,t)$ instead of $S(J,t)$ to emphasize that $\omega$
is varied now.
\end{remark}

\noindent
{\it Proof of Remark 3.3}\ \ Let $\omega + \delta\omega$ be a variation of the
K\"ahler form in a fixed K\"ahler class. Then $\delta\omega = \sqrt{-1}\partial
\barpartial \varphi$ for some real smooth function $\varphi$. By (19) the variation
$\delta S(\omega,t)$ of the perturbed scalar curvature is given by
\begin{eqnarray*}
&&\frac{\delta S(\omega,t)}{2m\pi} \omega^m 
+\frac{S(\omega,t)}{2m\pi}\,\Delta \varphi\,\omega^m \\
&=& \frac 1t (mc_m(\sqrt{-1}\partial\barpartial\varphi\otimes I + \frac {\sqrt{-1}}{2\pi}t\delta\Theta, 
\omega\otimes I + \frac {\sqrt{-1}}{2\pi}\,t\Theta, \\
&& \hspace{5cm}\cdots ,\omega\otimes I + \frac {\sqrt{-1}}{2\pi}\,t\Theta) 
- \Delta \varphi\,\omega^m)\\
&=& mc_m(\frac {\sqrt{-1}}{2\pi}\delta\Theta, 
\omega\otimes I + \frac {\sqrt{-1}}{2\pi}\,t\Theta, \cdots ,\omega\otimes I + \frac {\sqrt{-1}}{2\pi}\,t\Theta)\\
&& +\ m\,c_m(\sqrt{-1}\partial\barpartial\varphi\otimes I, 
\frac {\sqrt{-1}}{2\pi}\,\Theta,\omega\otimes I + \frac {\sqrt{-1}}{2\pi}\,t\Theta, \cdots ,\omega\otimes I + \frac {\sqrt{-1}}{2\pi}\,t\Theta)\\
&& + \cdots +  \ m\,c_m(\sqrt{-1}\partial\barpartial\varphi\otimes I, 
\frac {\sqrt{-1}}{2\pi}\,\Theta,\omega\otimes I , \cdots ,\omega\otimes I ).
\end{eqnarray*}
Thus
\begin{eqnarray*}
&&\frac 1{2m\pi}\delta(S(\omega,t)^2\,\omega^m) = 2S(\omega,t)mc_m(\frac {\sqrt{-1}}{2\pi}\delta\Theta, 
\omega\otimes I + \frac {\sqrt{-1}}{2\pi}\,t\Theta, \cdots ,\omega\otimes I + \frac {\sqrt{-1}}{2\pi}\,t\Theta)\\
&& +\ 2S(\omega,t)m\,c_m(\sqrt{-1}\partial\barpartial\varphi\otimes I, 
\frac {\sqrt{-1}}{2\pi}\,\Theta,\omega\otimes I + \frac {\sqrt{-1}}{2\pi}\,t\Theta, \cdots ,\omega\otimes I + \frac {\sqrt{-1}}{2\pi}\,t\Theta)\\
&& + \cdots +  \ 2S(\omega,t) m\,c_m(\sqrt{-1}\partial\barpartial\varphi\otimes I, 
\frac {\sqrt{-1}}{2\pi}\,\Theta,\omega\otimes I , \cdots ,\omega\otimes I ) \nonumber \\
&& - \ \frac 1{2m\pi}\ S(\omega,t)^2 \Delta\varphi\,\omega^m.
\end{eqnarray*}
Since $\delta\Theta = \nabla^{\prime\prime}\nabla^{\prime}\left(\varphi^i{}_j\right)$ we have
\begin{eqnarray}
&&\frac 1{2m\pi}\ \delta \int_M S(\omega,t)^2 \omega^m \\
&=& 2\int_M S(\omega,t)mc_m(\nabla_{\barl}\nabla_k \left(\varphi^i{}_j\right)\frac {\sqrt{-1}}{2\pi} 
d\overline{z^{\ell}}\wedge
dz^k, 
\omega\otimes I + \frac {\sqrt{-1}}{2\pi}\,t\Theta,\nonumber \\
&& \hspace{5cm} \cdots ,\omega\otimes I + \frac {\sqrt{-1}}{2\pi}\,t\Theta)\nonumber \\
&& +\ 2\int_M S(\omega,t)\ mc_m(\sqrt{-1}\partial\barpartial\varphi\otimes I, 
\frac {\sqrt{-1}}{2\pi}\,\Theta, \omega\otimes I + \frac {\sqrt{-1}}{2\pi}\,t\Theta, \nonumber\\
&& \hspace{5cm} \cdots ,\omega\otimes I + \frac {\sqrt{-1}}{2\pi}\,t\Theta)\nonumber \\
&& + \cdots + 2\int_M S(\omega,t)\ mc_m(\sqrt{-1}\partial\barpartial\varphi\otimes I, 
\frac {\sqrt{-1}}{2\pi}\,\Theta, \omega\otimes I, \cdots ,\omega\otimes I)\nonumber \\
&& -\ \frac 1{2m\pi}\ \int_m S(\omega,t)^2 \Delta\varphi\,\omega^m \nonumber
\end{eqnarray}
But
\begin{eqnarray}
\nabla_{\barl}\nabla_k\varphi^i{}_j &=& \nabla_{\barl}\nabla_k\nabla^i\,\varphi_j \\
&=& \nabla_{\barl}\nabla^i\nabla_k\varphi_j - \nabla_{\barl}(R^p{}_{jk}{}^i\,\varphi_p) \nonumber \\
&=& \nabla_{\barl}\nabla^i\nabla_k\varphi_j - (\nabla_{\barl}R^p{}_{jk}^{}i)\,\varphi_p 
- R^p{}_{jk}^{}i\,\varphi_{p\barl}  \nonumber \\
&=&  \nabla_{\barl}\nabla^i\nabla_k\varphi_j - (\nabla^pR_{\barl}{}_{jk}^{}i)\,\varphi_p 
- R^p{}_{jk}^{}i\,\varphi_{p\barl} \nonumber
\end{eqnarray}
where we used the second Bianchi identity at the last equality. It follows from (20) and (21) that
\begin{eqnarray}
&&\frac 1{2m\pi}\ \delta \int_M S(\omega,t)^2 \omega^m \\
&=& -2\int_M S(\omega,t)mc_m((\nabla_{\barl}\nabla^i\nabla_k \varphi_j
- \varphi_p\nabla^pR_{\barl jk}{}^i - R^p{}_{jk}{}^i\varphi_{p\barl})
\frac {\sqrt{-1}}{2\pi} 
dz^k \wedge d\overline{z^{\ell}}, \nonumber \\
&& \hspace{5cm}\omega\otimes I + \frac {\sqrt{-1}}{2\pi}\,t\Theta,
\cdots ,\omega\otimes I + \frac {\sqrt{-1}}{2\pi}\,t\Theta)\nonumber \\
&& +\ 2\int_M S(\omega,t)\ mc_m(\sqrt{-1}\partial\barpartial\varphi\otimes I, 
\frac {\sqrt{-1}}{2\pi}\,\Theta, \nonumber \\
&& \hspace{5cm}\omega\otimes I + \frac {\sqrt{-1}}{2\pi}\,t\Theta,
\cdots ,\omega\otimes I + \frac {\sqrt{-1}}{2\pi}\,t\Theta)\nonumber \\
&& \qquad+ \cdots + 2\int_M S(\omega,t)\ mc_m(\sqrt{-1}\partial\barpartial\varphi\otimes I, 
\frac {\sqrt{-1}}{2\pi}\,\Theta, \omega\otimes I, \cdots ,\omega\otimes I)\nonumber \\
&& -\ \frac 1{2m\pi}\ \int_m S(\omega,t)^2 \Delta\varphi\,\omega^m \nonumber
\end{eqnarray}
But 
$$ R_{\barl jk}{}^i \frac {\sqrt{-1}}{2\pi} dz^k\wedge d\overline{z^{\ell}} = 
R_{k\barl}{}^i{}_j \frac {\sqrt{-1}}{2\pi} dz^k\wedge d\overline{z^{\ell}} = \frac {\sqrt{-1}}{2\pi} \Theta.
$$
From this and integration by parts
\begin{eqnarray}
&&2 \int_M S(\omega,t)\,m\,c_m(\varphi_p\nabla^pR_{\barl jk}{}^i \frac {\sqrt{-1}}{2\pi} 
dz^k \wedge d\overline{z^{\ell}}, \omega\otimes I + \frac {\sqrt{-1}}{2\pi}\,t\Theta,\\
&& \hspace {2cm}\cdots ,\omega\otimes I + \frac {\sqrt{-1}}{2\pi}\,t\Theta)
 = - \frac 1{2m\pi}\ \int_M S(\omega,t)^2 \Delta\varphi\ \omega^m. \nonumber
\end{eqnarray}
It follows from (22) and (23) that
\begin{eqnarray}
&&\ \frac 1{2m\pi}\ \delta \int_M S(\omega,t)^2 \omega^m \nonumber \\
&=& -2\int_M S(\omega,t)mc_m((\nabla_{\barl}\nabla^i\nabla_k \varphi_j
\frac {\sqrt{-1}}{2\pi} 
dz^k \wedge d\overline{z^{\ell}}, \nonumber \\
&& \hspace{5cm}\omega\otimes I + \frac {\sqrt{-1}}{2\pi}\,t\Theta,
\cdots ,\omega\otimes I + \frac {\sqrt{-1}}{2\pi}\,t\Theta)\nonumber \\
&& + 2\int_M S(\omega,t)mc_m(R^p{}_{jk}{}^i\varphi_{p\barl}
\frac {\sqrt{-1}}{2\pi} 
dz^k \wedge d\overline{z^{\ell}}, \\
&& \hspace{5cm}\omega\otimes I + \frac {\sqrt{-1}}{2\pi}\,t\Theta,
\cdots ,\omega\otimes I + \frac {\sqrt{-1}}{2\pi}\,t\Theta) \nonumber\\
&& +\ 2\int_M S(\omega,t)\  mc_m(\sqrt{-1}\partial\barpartial\varphi\otimes I, 
\frac {\sqrt{-1}}{2\pi}\,\Theta, \\
&& \hspace{5cm}\omega\otimes I + \frac {\sqrt{-1}}{2\pi}\,t\Theta,
 \cdots ,\omega\otimes I + \frac {\sqrt{-1}}{2\pi}\,t\Theta)\nonumber  \\
 && \qquad\qquad + \cdots + 2\int_M S(\omega,t)\  mc_m(\sqrt{-1}\partial\barpartial\varphi\otimes I, 
\frac {\sqrt{-1}}{2\pi}\,\Theta, \omega\otimes I ,
 \cdots ,\omega\otimes I ) \nonumber \\
&& -\ \frac 1{m\pi}\ \int_m S(\omega,t)^2 \Delta\varphi\,\omega^m 
\end{eqnarray}
When $t=0$ this is equal to
\begin{eqnarray}
\ \frac 1{2m\pi}\ \delta \int_M S^2 \omega^m 
&=& -2\int_M S\,{\overline D}\varphi\,\omega^m
 + 2\int_M S\,\sum_{i,j=1}^m \ \frac 1{2\pi}R_{i\barj}\varphi^{i\barj} \omega^m
 \\
&& +\ 2\int_M S\  \sum_{i\ne j}\varphi_{i\bari}\ \frac 1{2\pi}\,R_{j\barj}\,\omega^m 
 -\ \frac 1{m\pi}\ \int_m S^2 \Delta\varphi\,\omega^m \nonumber
\end{eqnarray}
with $D = \nabla_i\nabla_j\nabla^i\nabla^j$ where $S = S(\omega,0)$ is the unperturbed scalar curvature and we used the normal
coordinates such that the complex Hessian $(\varphi_{i\barj})$ is diagonalized. 
The third term on the right hand side can then be computed using
\begin{eqnarray*}
\sum_{i\ne j} \varphi_{i\bari}\ \frac 1{2\pi}\,R_{j\barj} &=& (\sum_{i=1}^m \varphi_{i\bari})(\sum_{j=1}^m\ \frac 1{2\pi}\, R_{j\barj}) - 
\varphi^{i\barj}\ \frac 1{2\pi}\, R_{i\barj} \\
&=& \Delta \varphi\ \frac 1{2m\pi}\,S - \varphi^{i\barj}\ \frac 1{2\pi}\, R_{i\barj},
\end{eqnarray*}
and we see from this and (27) that 
$$
\ \frac 1{2m\pi}\,\delta \int_M S^2 \omega^m 
= -2\int_M DS\,\varphi\,\omega^m.
$$
This proves the fact that the critical points in the unperturbed case are the extremal K\"ahler metrics.
We have seen that when $t=0$, (24) + (25) + (26) vanishes. But when $t \ne 0$, this is not the case
because we have the term with $t^{m-1}$ only in (24)
$$2\int_M S(\omega,t)mc_m(R^p{}_{jk}{}^i\varphi_{p\barl}
\frac {\sqrt{-1}}{2\pi} 
dz^k \wedge d\overline{z^{\ell}}, \frac {\sqrt{-1}}{2\pi}\,t\Theta,
\cdots ,\frac {\sqrt{-1}}{2\pi}\,t\Theta),$$
which does not always vanish. This completes the proof of Remark 3.3.

\section{K\"ahler metrics of harmonic Chern forms}

Let $M$ be a compact K\"ahler manifold with a fixed K\"ahler class $[\omega_0]$
and $\frak h(M)$ the complex Lie algebra of all holomorphic vector fields.
For any $\omega \in [\omega_0]$, let $c_k(\omega)$ be the $k$-th Chern form 
with respect to $\omega$ as in Remark 3.3. 
Let 
$Hc_k(\omega)$ be the harmonic part of $c_k(\omega)$. Here the harmonic projection $H$ 
is taken with respect to the K\"ahler metric $\omega$. Then
$$ c_k(\omega) - Hc_k(\omega) = \sqrt{-1} \partial \barpartial F_k $$
for some smooth real $(k-1,k-1)$-form
$$ F_k \in \Omega^{k-1,k-1}(M). $$
For a holomorphic vector field $X \in {\frak h}(M)$, 
define $f_k : \frak h(M) \to \bfC$ by
$$ f_k(X) = \frac 1{m-k+1}\int_M L_XF_k \wedge \omega^{m-k+1}.$$
\begin{theorem}[S. Bando \cite{bando83}]\ \ The functional $f_k$ on 
${\frak h}(M)$ is independent of the choice of
$\omega \in [\omega_0]$, becomes a Lie algebra character 
and obstructs the existence of K\"ahler metrics
$\omega$ in $[\omega_0]$ of 
harmonic $k$-th Chern form.
\end{theorem}
\noindent
In \cite{futaki04} the author gave a larger family of integral invariants including
$f_i$'s and obstructions to asymptotic Chow semi-stability.

Here again as in Remark 3.3 we are fixing $J$ and varying $\omega$, 
instead of fixing $\omega$ and varying $J$.
So we denote the perturbed scalar curvature by $S(\omega,t)$ as in (19).
If $X = \mathrm{grad}^{\prime} u = 
g^{i\barj}\frac{\partial u}{\partial \overline{z^j}}\frac{\partial}{\partial z^i}$ with 
$\int_M u\ \omega^m = 0$ then we see using the integration by parts that
\begin{equation}
\ \frac 1{2m\pi}\,\int_M u\,S(\omega,t)\ \omega^m = -\ f_1(X)\ -\ t\,f_2(X)
 -\ \cdots\ -\ t^{m-1}f_m(X).
\end{equation}
We put
$$ F_t(X) := f_1(X) + t\,f_2(X) + \cdots + t^{m-1}f_m(X).$$
and call it {\bf total Bando character}.\\

\medskip

\noindent
\begin{proposition}\ \ For fixed small  $t \in \bfR$, $F_t : \frak h(M) \to \bfC$ is an obstruction to
to the existence of K\"ahler metric  $\omega \in [\omega_0]$ 
of constant perturbed scalar curvature $S(\omega, t)$. If there exists a perturbed extremal K\"ahler
metric and the total Bando character vanishes, then the perturbed extremal K\"ahler
metric has constant perturbed scalar curvature.
\end{proposition}

\begin{proof} If there is a K\"ahler form $\omega \in [\omega_0]$ such that
$S(\omega,t)$ is constant. Then the total Bando character has to vanish because
of (28) and the normalization $\int_M u\ \omega^m = 0$. If $\omega$ is a 
perturbed extremal metric then $\mathrm{grad}' S(\omega,t)$ is a holomorphic vector field and
$$
F_t(\mathrm{grad}' S(\omega,t)) = \ \frac 1{2m\pi}\,\int_M g^{i\barj}\frac{\partial S(\omega,t)}{\partial z^i}
\frac{\partial S(\omega,t)}{\partial \barz^j}\ \omega^m.
$$
Thus if $F_t$ vanishes then $S(\omega,t)$ is constant.
\end{proof}

Let $\sigma(t)$ be the topological invariant
$$ \sigma(t) = \frac {(c_1(M)\wedge [\omega_0]^{m-1} + t\,c_2(M)[\omega_0]^{m-2} + 
\cdots + t^{m-1}c_m(M))[M]}{[\omega_0]^m[M]}.$$
This is obviously the average of the perturbed scalar curvature (with respect to any K\"ahler
form $\omega \in [\omega_0]$). For any two K\"ahler forms $\omega^{\prime}$ and
$\omega^{\prime\prime}$ we define
$$ \mathcal M_t(\omega^{\prime},\omega^{\prime\prime}) =
- \int_0^1 ds\int_M \frac{\partial \varphi_s}{\partial s} (S(\omega_s,t) - \sigma(t))\omega_s^m$$
where $\omega_s = \omega + \sqrt{-1}\partial\barpartial\,\varphi_s$, $0\le s \le 1$, 
is a smooth path in $[\omega_0]$ joining
$\omega^{\prime}$ and $\omega^{\prime\prime}$. Bando and Mabuchi (\cite{bandomabuchi86})
observed that every coefficient of $t^k$ in $\mathcal M_t(\omega^{\prime},\omega^{\prime\prime})$,
and thus $\mathcal M_t(\omega^{\prime},\omega^{\prime\prime})$, is independent of the
choice of the paths $\omega_s$ and satisfies the cocycle conditions.
Putting $\nu_t(\omega) := \mathcal M_t(\omega_0,\omega)$, we get a functional on the 
space of all K\"ahler forms in the cohomology class $[\omega_0]$. The functional $\nu_0$ 
in the case when $t=0$ is the so-called K-energy or Mabuchi energy.
We call $\nu_t$ the perturbed Mabuchi energy. It is obvious that the critical points of the
perturbed Mabuchi energy are the K\"ahler metrics of constant perturbed scalar curvature.
In the case when $t=0$ Chen and Tian \cite{chentian04} proved that the Mabuchi energy is bounded from
below if there exists a K\"ahler metric of constant scalar curvature, and that the infimum of
the Mabuchi energy is attained exactly on the space of K\"ahler metrics of constant
scalar curvature, extending earlier result of Bando and Mabuchi \cite{bandomabuchi87}
for K\"ahler-Einstein manifolds of positive first Chern class. We hope to discuss 
for the perturbed case in a later paper.

The proof of the fact that the definition of $\mathcal M_t$
is independent of the paths follows from the fact that $S(\omega,t)\omega^m$ gives a
closed 1-form on the space of K\"ahler forms. The closedness comes from the symmetry
between the endomorphism part and the form part in the definition of $S(\omega,t)\omega^m$,
as was explained between the equation (9) and Theorem 2.2. The detailed discussion was given in \cite{futaki88}
but of course the original idea goes back to Bando \cite{bando83}.

For the identity component $\mathrm{Aut}^0(M)$ of the group of all holomorphic automorphisms of $M$,
let $G$ denote the maximal linear algebraic subgroup. The maximal reductive subgroup $K^c$ of $G$ is
the complexification of a compact Lie group $K$. Taking the average of the K\"ahler metric by the action 
of $K$ we may assume that $K$ acts as isometries. We denote by $\omega$ the K\"ahler form of the
averaged K\"ahler metric. Then the elements of the Lie algebra of $K$ are Killing vector fields of
$(M, \omega)$ 
and are thus obtained as the real parts of the gradient vector fields of purely imaginary functions (see e.g.
\cite{kobayashitrnsf}). Therefore as a complex Lie algebra, the Lie algebra $\frak k^c$ is isomorphic
to the Lie algebra $\frak u$ spanned over $\bfC$ by some 
real functions $u_1, \cdots, u_d$ with the normalization $\int_Mu_i\,\omega^m = 0$ where the Lie bracket on $\frak u$ 
is given by the Poisson bracket
$$ \{u,v\} = u^iv_i - v^iu_i = g^{i\barj}\frac{\partial u}{\partial \barz^j}\frac{\partial v}{\partial z^i}  - 
g^{i\barj}\frac{\partial v}{\partial \barz^j} \frac{\partial u}{\partial z^i}.$$
\begin{proposition}\ \ Let the situation be as above. 
If we choose $\omega_r = \omega + \sqrt{-1}\partial\barpartial\,\varphi_r$
so that $\varphi_0 = 0$ and that $\dot{\varphi}_r|_{r=0} = u$ for some
real smooth function $u$ in $\frak u$, then
$$\left. \frac d{dr}\right|_{r=0}\nu_t(\omega_r) = 2m\pi\,F_t(\mathrm{grad}^{\prime}u).$$
\end{proposition}
\begin{proof}\ \ This is immediate from
\begin{equation*}
\nu_t(\omega_r) = - \int_0^r dq \int_M \frac{\partial\varphi_q}{\partial q} (S(\omega_q,t) -
\sigma(t))\omega_q^m
\end{equation*}
and
\begin{eqnarray*}
\left. \frac d{dr}\right|_{r=0}\nu_t(\omega_r) &=& 
- \int_M u\ S(\omega, t)\ \omega^m \\
&=& 2m\pi\,F_t(\mathrm{grad}^{\prime}u)
\end{eqnarray*}
where the last equality follows because $u$ is a normalized Hamiltonian function for
a holomorphic vector field.
\end{proof}
\noindent
This proposition shows that the perturbed Mabuchi energy is an integral form of the total
Bando character. A way of computing the unperturbed Mabuchi energy $\nu_0$ without using the path integral
was given in \cite{futakinakagawa}. It would be interesting if one can give a formula for $\nu_t$
without using path integral. B. Weinkove \cite{weinkove02} related the degree $1$ and $2$ terms in $t$ of
$\mathcal M_t$ to
Donladson's functional which was used in the proof of the existence of Hermitian-Einstein metrics on
stable vector bundles \cite{donaldson85}.

We also remark that the modified Mabuchi energy to treat the extremal metrics can be also defined
in the perturbed case just as defined in \cite{guanchen} and \cite{simanca00}. One can use the
proof given in \cite{guan99}.

The results obtained above 
may be interesting to compare with a results of X. Wang \cite{xwang04} (see also \cite{futaki05})
which we summarize below.

Let $(Z, \Omega)$ be a K\"ahler manifold and 
suppose a compact Lie group $K$ acts on $Z$ as holomorphic isometries. 
Then the complexification $K^c$ of $K$ also acts on $Z$ as biholomorphisms. 
The actions of $K$ and $K^c$ induce homomorphisms of the 
Lie algebras ${\frak k}$ and ${\frak k}^c$ to the real Lie algebra $\Gamma(TZ)$ of all 
smooth vector fields on $Z$, both of which we denote by $\rho$. 
If $\xi + i\eta \in {\frak k}^c$ with $\xi,\ \eta \in {\frak k}$, then 
$$ \rho(\xi + i\eta) = \rho(\xi) + J\rho(\eta),$$
where $J$ is the complex structure of $Z$.
Suppose $[\Omega]$ is an integral class and there is a holomorphic line 
bundle $L \to Z$ with $c_1(L) = [\Omega]$. 
There is an Hermitian metric $h$ of $L^{-1}$ such that 
its Hermitian connection $\theta$ satisfies 
$$ - \frac 1{2\pi} d\theta = \Omega.$$
Suppose we have a lifting of $K^c$ to $L^{-1}$, so that we have a moment map
$\mu : Z \to {\frak k}^{\ast}$ because the lifting of $K$-action to $L$ is equivalent
to defining a moment map (see \cite{donkro}, section 6.5). 
Let $\pi : L^{-1} \to Z$ be 
the projection and $\pi (p) = x$ with $p \in L^{-1} - \mathrm{zero\ section},\ x \in Z$. 
Denote by $\Gamma = K^c\cdot x$ the $K^c$-orbit of $x$ in $Z$, 
and $\widetilde{\Gamma} = K^c\cdot p$ be the $K^c$-orbit of $p$ in $L^{-1}$.
We say that $x \in Z$ is polystable with respect to the $K^c$-action if the orbit $\widetilde{\Gamma}$ is closed in $L^{-1}$.
Consider the function $h : \widetilde{\Gamma} \to {\mathbb R}$ defined by 
$$ h(\gamma) = \log |\gamma|^2. $$
Fundamental facts are
\begin{itemize}
\item\ \ $h$ has a critical point if and only if 
the moment map $\mu : Z \to \frak k^{\ast}$ has a zero along $\Gamma$:
\item\ \ $h$ is a convex function.
\end{itemize}
For these facts refer again to \cite{donkro}, section 6.5. These imply the following two 
propositions.

\begin{proposition}
A point $x \in Z$ is polystable with respect to the action of $K^c$ if and only if 
the moment map $\mu$ has a zero along $\Gamma$.
\end{proposition}

\begin{proposition}
The set $\{x \in \Gamma\ |\ \mu(x) = 0\}$ has only one component, and the 
orbit $Stab(x)^c\cdot x$ of the complexification of the stabilizer at $x$ through $x$ is 
connected even if $Stab(x)^c$ is not connected.
\end{proposition}

For a given $x \in Z$ we extend $\mu(x) : {\frak k} \to \bfR$ complex linearly
to $\mu(x) : {\frak k}^c \to \bfC$. 
For notational convenience we denote by $K_x$ (resp. $(K^c)_x$) 
the stabilizer of $x$ in $K$ (resp. $K^c$), and by ${\frak k}_x$ 
and $({\frak k}^c)_x$ the Lie algebra of $K_x$ and $(K^c)_x$. 
Define $f_x : ({\frak k}^c)_x \to {\mathbb C}$ to be the restriction of 
$\mu(x) : {\frak k}^c \to \bfC$ to $({\frak k}^c)_x$. 
Note that $(K^c)_{gx} = g(K^c)_x g^{-1}$. 

\begin{proposition}[Wang \cite{xwang04}]
Fix $x_0 \in Z$. Then for $x \in K^c\cdot x_0$, 
$f_x$ is $K^c$-equivariant in that $f_{gx}(Y) = f_x(Ad(g^{-1})Y)$. In particular 
if  $f_x$ vanishes at some $x \in K^c\cdot x_0$ it vanishes at
all $x \in K^c\cdot x_0$. 
Moreover $f_x : ({\frak k}^c)_x \to {\mathbb C}$ is a Lie algebra character.
\end{proposition}

For a proof of this proposition, see \cite{xwang04} and also \cite{futaki05}. 
Suppose now we are given a $K$-invariant inner product on ${\frak k}$. 
Then we can identify ${\frak k} \cong {\frak k}^{\ast}$, and ${\frak k}^{\ast}$ 
has a $K$-invariant inner product. Consider the function $\phi : K^c\cdot x_0 \to {\mathbb R}$ 
defined by
$\phi(x) = |\mu(x)|^2$. We say that $x \in K^c\cdot x_0$ is an extremal point if $x$ is 
a critical point of $\phi$.  

\begin{proposition}[Wang \cite{xwang04}]
Let $x \in K^c\cdot x_0$ be an extremal point. 
Then we have a decomposition
$$ ({\frak k}^c)_x = ({\frak k}_x)^c \oplus \sum_{\lambda > 0} {\frak k}^c_{\lambda} $$
where ${\frak k}^c_{\lambda}$ is $\lambda$-eigenspace of ${\mathrm ad}(i\mu(x))$, 
and $i\mu(x)$ lies in the center of $({\frak k}_x)^c$. In particular 
$({\frak k}_x)^c = ({\frak k}^c)_x$ if and only if $\mu(x) = 0$.
\end{proposition}

For a proof of this proposition, see \cite{xwang04} and also \cite{futaki05}. 
Let $(M, \omega_0, J_0)$ be a compact K\"ahler manifold with a fixed K\"ahler form $\omega_0$. 
Apply the above results for  finite dimensional manifold $Z$ to the set $\mathcal J$ of all $\omega$-compatible integral complex structures $J$ with 
respect to which $(M, \omega_0, J)$ is a K\"ahler manifold, where 
the compact Lie group $K$ is replaced by the group of 
symplectomorphisms generated by Hamiltonian diffeomorphisms. This explains a relationship
between stability and various results about extremal K\"ahler metrics. 
For example, Proposition 4.6 explains the total Bando character 
and Proposition 4.7 of course explains Calabi's decomposition
theorem for the Lie algebras of all holomorphic vector fields on compact extremal K\"ahler manifolds
\cite{calabi85} (see the next section).

\section{Deformations of extremal K\"ahler metrics}

Let $M$ be a compact complex manifold carrying a K\"ahler metric. By a
($t$-perturbed) extremal K\"ahler 
class we mean a de Rham cohomology class which contains the K\"ahler form of a
($t$-perturbed) extremal K\"ahler metric. In this section we prove the following
result which extends the results of LeBrun and Simanca \cite{lebrunsimanca93}, 
\cite{lebrunsimanca94}.
\begin{theorem}\ \ For an extremal K\"ahler class $[\omega_0]$, there
exists a neighborhood $U\times (-\epsilon,\epsilon)$ of 
$([\omega_0],t)$ in $H^{1,1}_{DR}(M,\bfR) \times \bfR$
such that all points of $U$ are $t$-perturbed extremal K\"ahler classes for
all $t \in (-\epsilon, \epsilon)$.
\end{theorem}

The rest of this section is devoted to the proof of this theorem. We first review
well known facts on Hamiltonian holomorphic vector fields on compact
K\"ahler manifolds. Let $(M,g)$ be a compact K\"ahler manifold. We
define a fourth-order elliptic differential operator $L_g : C^{\infty}_{\bfC}(M) \to
C^{\infty}_{\bfC}(M)$ by
$$ L_g\,u = \nabla^{\prime\prime\ast}\nabla^{\prime\prime\ast}
\nabla^{\prime\prime}\nabla^{\prime\prime}\,u, $$
where $C^{\infty}_{\bfC}(M)$ denotes the set of all complex valued smooth
functions on $M$. More precisely
\begin{eqnarray*}
 L_g\,u &=& \nabla^{\barj}\nabla^{\bari}\nabla_{\barj}\nabla_{\bari}\,u \\
 &=& \Delta^2u + R^{\barj i}\nabla_{\barj}\nabla_i\,u + \nabla^{\barj}S\,\nabla_{\barj}u
 \end{eqnarray*}
 where $S$ denotes the unperturbed scalar curvature.
 Then the kernel of $L_g$ consists of all smooth functions $u$ whose gradient vector fields
 $$ \mathrm{grad}^{\prime} u := g^{i\barj}\nabla_{\barj}u\,\frac{\partial}{\partial z^i}$$
 are holomorphic vector fields. It is well known that such holomorphic vector fields are
 exactly those which have zeros (see \cite{lebrunsimanca93} for a comprehensive proof).
 Since constant functions correspond to the zero vector field, we only consider the subspace
 $(\mathrm{ker}\, L_g)_0$ consisting of all functions $u \in \mathrm{ker}\, L_g$ which are
 orthogonal to constant functions:
 $$ \int_M u\,\omega_g^m = 0.$$
 Now we study the behavior of $u \in (\mathrm{ker} L_g)_0$ when the K\"ahler metric $g$
 varies in the same K\"ahler class. The following lemma was used in \cite{futakimabuchi95},
 pp.208-209, but we will reproduce a proof here for the reader's convenience.
 
 \begin{lemma}\ \ Let $\tildeg_{i\barj} = g_{i\barj} + \nabla_i\nabla_{\barj}\,\varphi$ be
 a K\"ahler metric in the same K\"ahler class as $g_{i\barj}$. If 
 $u \in (\mathrm{ker}\,L_{\tildeg})_0$, then $\tildeu := u + \nabla^iu\,\nabla_i\varphi \in
 (\mathrm{ker}\,L_{\tildeg})_0$ and $\mathrm{grad}_{\tildeg}\tildeu = \mathrm{grad}_g\,u$.
 \end{lemma}
 \begin{proof}\ \ We first show the last equation.
 \begin{eqnarray*}
 \mathrm{grad}_{\tildeg}\,\tildeu &=& \tildeg^{i\barj}\,\frac{\partial\tildeu}{\partial\barz^j}
\, \frac{\partial}{\partial z^i} = 
\tildeg^{i\barj}\,(\frac{\partial u}{\partial\barz^j} + \nabla^ku\nabla_k\nabla_{\barj}\varphi)
\, \frac{\partial}{\partial z^i}\\
&=& \tildeg^{i\barj}\,\nabla^ku\,
(g_{k\barj} + \nabla_k\nabla_{\barj}\varphi)
\, \frac{\partial}{\partial z^i} = \nabla^iu\,\frac{\partial}{\partial z^i} = \mathrm{grad}_gu.
\end{eqnarray*}
It remains to see
$$ \int_M \tildeu\,\omega_{\tildeg}^m = 0.$$
Let $g_{ti\barj} = g_{i\barj} + t\nabla_i\nabla_{\barj}\,\varphi$ be the line segment of
K\"ahler metrics between $g$ and $\tildeg$, and $u_t = u + t\nabla^iu\,\nabla_i\varphi$ be the corresponding functions in $(\mathrm{ker}\,L_{\tildeg})_0$. It is sufficient to prove
$$ \frac{d}{dt}\int_M\,u_t\,\omega_{g_t}^m = 0.$$
It is also sufficient to prove this at $t=0$. But
$$ \left.\frac{d}{dt}\right|_{t=0}\,\int_M u_t\,\omega_{g_t}^m = 
\int_M (\nabla^iu\nabla_i\varphi + u(\Delta\varphi))\,\omega_g^m = 0,$$
where $\Delta = \nabla^i\nabla_i\,u$ denotes the complex Laplacian. This 
completes the proof.
\end{proof}
Now let $K$ be the identity component of the isometry group of $(M,g)$, and
$\frak k$ be its Lie algebra. Hence $\frak k$ consists of all Killing vector fields.
On a compact K\"ahler manifold $\frak k$ can be embedded into the complex
Lie algebra ${\frak h}(M)$ of all holomorphic vector fields on $M$ by $X \in {\frak k} \mapsto
\frac12(X - \sqrt{-1}JX) \in {\frak h}(M)$. By this $\frak k$ is often identified with the 
image in ${\frak h}(M)$ of this embedding. As was explained in the previous section
when a holomorphic vector field $X$ is written as a gradient vector field of a
complex valued smooth function, $X$  is a Killing vector field if and only if the 
function is a purely imaginary valued function. We choose real valued smooth functions
$u_1,\ \cdots,\ u_d$ so that the gradient vector fields of $iu_1,\ \cdots,\ iu_d$ form a basis
of $\frak k \otimes\bfC$. We also assume that $1,\ u_1,\ \cdots,\ u_d$ form an 
$L^2$-orthonormal system (under the normalization $\int_M u_i \omega^m = 0$).
Let us denote by $J_g$ the linear span over $\bfC$ of  $1,\ u_1,\ \cdots,\ u_d$. 
\begin{remark}\ \ Since the imaginary part of $\mathrm{grad}^{\prime}\,u_j$ is a Killing vector field,
$(\mathrm{grad}^{\prime}\,u_j)\varphi$ is a real function for a $K$-invariant real function $\varphi$.
\end{remark}
\begin{remark}\ \ If $\tildeg_{i\barj} = g_{i\barj} + \nabla_i\nabla_{\barj}\,\varphi$ is a
$K$-invariant K\"ahler metric in the same K\"ahler class as $g$, then the corresponding basis of 
$J_{\tildeg}$ consisting of real functions are 
$$ 1,\ \tildeu_1= u_1 + (\mathrm{grad}^{\prime}\,u_1)\varphi.,\ \cdots\ ,
\ \tildeu_d = u_d + (\mathrm{grad}^{\prime}\,u_d)\varphi.$$
It is easy to see that they form an $L^2$-orthonormal system with respect to $\tildeg$ (see
\cite{futakimabuchi95}, Appendix 2).
\end{remark}

Since we assume that there is an extremal K\"ahler metric, the Lie algebra ${\frak h}(M)$
has the following structure by a theorem of Calabi \cite{calabi85}. Namely there is a
decomposition
$$ {\frak h}(M) = {\frak h}_0 + \sum_{\lambda \ne 0} {\frak h}_{\lambda}, $$
where ${\frak h}_{\lambda}$ is a $\lambda$-eigenspace of the adjoint action
of the extremal vector field 
$$\mathrm{ad}(
\mathrm{grad}^{\prime}\,S) : {\frak h}(M) \to {\frak h}(M),$$
and further ${\frak h}_0$ 
is the complexification of the Lie algebra $\frak k$ consisting of all Killing vector fields
on $(M,g)$. In particular, it turns out that $\mathrm{grad}^{\prime}\,S$ lies in the center
of ${\frak h}_0$. that $[{\frak h}_{\lambda},{\frak h}_{\mu}] \subset {\frak h}_{\lambda + \mu}$,
and that ${\frak h}_0$ is a maximal reductive Lie subalgebra of ${\frak h}(M)$.

Now we consider the set of all K\"ahler metrics invariant under the identity component
of the isometry group $K$ of $(M,g)$ of the form
$$ \omega(\alpha, \varphi) = \omega + \alpha + \sqrt{-1}\partial\barpartial\varphi$$
where $\alpha$ is a $K$-invariant real harmonic $(1,1)$-form on $(M,g)$ and $\varphi$
is a $K$-invariant real-valued $L^2_{k+4}$-function. Hence the space of such $K$-invariant
K\"ahler metrics is identified with an open subset of $H^{1,1}(M;\bfR)\, \times\,
L^2_{k+4,K}$ where $H^{1,1}(M;\bfR)$ denotes the vector space of all
real harmonic $(1,1)$-forms on $M$ and $L^2_{k+4,K}$ is the vector space of all real
valued $K$-invariant $L^2_{k+4}$ functions on $M$. Let $I_{k+4}$ be the orthogonal
complement to the subspace spanned by $1,\ u_1,\ \cdots,\ u_d$ in $L^2_{k+4,K}$.

Let $\tildeg$ be the K\"ahler metric corresponding to $\omega(\alpha, \varphi).$
Then we obtain, as in Remark 5.4, $L^2_{k+3}$-functions $(1,\,\tildeu_1,\,\cdots,\,\tildeu_d)$
whose gradient vector fields span the Lie algebra ${\frak k}$. Let $\widetilde{J}_{k+3}$ be the
linear span of $(1,\,\tildeu_1,\,\cdots,\,\tildeu_d)$. We put $\tildeu_0 = 1$.
Then for a sufficiently small neighborhood $U$ of $g$ in $H^{1,1}(M;\bfR) \times L^2_{k+4,K}$,
we have 
$$\det(u_i,\tildeu_j)_{L^2} \ne 0$$
for all $\tildeg \in U$. Then it is easy to see
$$ \mathrm{ker}(1 - \Pi_g)(1 - \Pi_{\tildeg}) = \mathrm{ker}(1 - \Pi_{\tildeg})$$
where $\Pi_g$ and $\Pi_{\tildeg}$ are respectively the $L^2$ projections of 
$L^2_{k,K}$ onto $J_{k+3} \subset L^2_{k,K}$ and onto $\widetilde{J}_{k+3} \subset L^2_{k,K}$:
$$ \Pi_g : L^2_{k,K} \to L^2_{k,K},\qquad\qquad \Pi_g(f) = \sum_{i=0}^d (f,u_i)u_i, $$
$$ \Pi_{\tildeg} : L^2_{k,K} \to L^2_{k,K},\qquad\qquad 
\Pi_{\tildeg}(f) = \sum_{i=0}^d (f,\tildeu_i)\tildeu_i.$$

Put $V := U \cap (H^{1,1}(M;\bfR)\times I_{k+4})$, and take a neighborhood $W$ of the origin
in $V \times \bfR$ such that for every point $(\tildeg,t)$ in $W$ (identifying $V$
with the space of K\"ahler metrics)
the inner product (9) makes sense so that one can consider $t$-perturbed scalar curvature.
Consider the map $\frak S : W \to I_k$
defined by
$$ {\frak S}(\tildeg,t) = (1 - \Pi_g)(1 - \Pi_{\tildeg})S(\tildeg,t).$$
Note that ${\frak S}(g,0) = 0$ and that ${\frak S}^{-1}(0)$ is the set of all perturbed
extremal K\"ahler metrics in $W$. To complete the proof of Theorem 5.1, it is sufficient to
show, by the implicit function theorem, that the partial derivative
$$ D{\frak S}_{(g,0)} : I_{k+4} \to I_k $$
at $(g,0)$ in the direction of $I_{k+4}$ is an isomorphism. In the direction of $\psi \in I_{k+4}$,
the derivative of the scalar curvature is
$$ (DS)_{g}(\psi) = - \Delta^2\psi - R^{\barj i}\nabla_{\barj}\nabla_i\,\psi,$$
and the derivative of the projection $\Pi$ is
\begin{eqnarray*}
(D\Pi)(S(g))_g(\psi) &=& \left.\frac d{dt}\right|_{t=0} (S + \nabla^iS\,t\nabla_i\psi)\\
&=& \nabla^iS\nabla_i\psi = \nabla^{\bari}S\nabla_{\bari}\psi,
\end{eqnarray*}
where the last equality follows from Remark 5.3. Combining these two equations, we obtain
\begin{eqnarray*}
(D{\frak S})_g (\psi) &=& (1-\Pi_g)(-\Delta^2\psi - R^{\barj i}\nabla_{\barj}\nabla_i\psi - 
\nabla^{\barj}S\nabla_{\barj}\psi)\\
&=& (1-\Pi_g)(-L_g\psi)
\end{eqnarray*}
If $(1-\Pi_g)(L_g\psi)=0$, then $L_g\psi \in J_g$. But since $L_g$ is self-adjoint,
$(\mathrm{Image}\,L_g)^{\perp} = \mathrm{ker}\,L_g$ and hence $L_g \psi = 0$.
Since $\psi \in I_{k+4}$, this implies $\psi = 0$. Thus $(D{\frak S})_{(g,0)}$ is injective,
which also implies that $(D{\frak S})_{(g,0)}$ is surjective since $(D{\frak S})_{(g,0)}$ is
self-adjoint. This completes the proof.

\end{document}